\documentclass[12pt]{article}

\usepackage{amsmath}
\usepackage{amssymb}
\usepackage{amsthm}
\usepackage{newtxtext,newtxmath}
\usepackage{url}

\theoremstyle{definition}
\newtheorem{theorem}{Theorem}[section]

\title{On Criteria about Estimation of the Riemann Zeta Function on the Line $\sigma=1$}
\date{5 May, 2020}
\author{Yoshihiro Koya}

\begin{document}
\maketitle

\begin{abstract}
  In this paper we give criteria about estimation of derivatives of the
  Riemann Zeta Function on the line $\sigma=1$.
\end{abstract}

\section{Introduction}

Let $s = \sigma + i t$.
The Riemann zeta function is defined by the Dirichlet series
\begin{equation*}
  \zeta(s) = \sum^{\infty}_{n=1}\frac{1}{n^s} \qquad (\sigma > 1).
\end{equation*}
This series is absolutely convergent and analytic for $\sigma > 1$.
$\zeta(s)$ is analytically continued to the whole plane, and regular for all $s$ except for a simple pole at $s=1$.

It is known that $\zeta(s)$ has zeros at $s = -2, -4, \ldots$, which are called \textit{trivial zeros}.
As for the zeros of $\zeta(s)$,
there is the celebrated conjecture, \textit{the Riemann hypothesis}.
It states that all zeros of $\zeta(s)$ other than the trivial zeros are on
the line $\sigma =\frac{1}{2}$.

There are many facts which are derived from the Riemann hypothesis.
Among such, there are also some conjectures about estimations of the values of $\zeta(s)$ itself.
For example, we can expect the following big-$O$ results~:
\begin{align*}
  \zeta(1 + it) &= O(\log\log t) \\
  \noalign{and}
  \frac{\zeta^\prime(1+it)}{\zeta(1+it)} &= O(1),
\end{align*}
by assuming the Riemann hypothesis(\cite[Theorem~14.8 and Theorem~14.5]{Titch}).
Note that it seems to be quite hard to obtain these results without assuming
the Riemann hypothesis.

\bigskip
In this paper, we shall give some criteria which give
the relationship beetween 
the estimations of the values and those of some integrals of $\zeta^{(m)}(1+it)$.
More exactly, we shall prove~:

\begin{theorem}[cf.~Theorem~13.2 of \cite{Titch}] \label{thm:main}
  Let $\epsilon > 0 $ be an arbitrary positive number.
  For a fixed non-negative integer $m$, the followings are equivalent~:
  \begin{enumerate}
  \item $\zeta^{(m)} (1 + it) = O(\log^\epsilon t)$.
  \item For each positive integer $k = 1, 2, \ldots$,
    \begin{equation*}
      \dfrac{1}{2\Delta} \int^{T+\Delta}_{T-\Delta} |\zeta^{(m)}(1+it)|^{2k} dt = O(\log^\epsilon T),
    \end{equation*}
    where $\Delta = 1/\log^{m+2} T$.
  \end{enumerate}
\end{theorem}

\section{Criteria about Estimation of $\zeta^{(m)}(1+it)$}

In this paper, we shall prove the following theorem.

\begin{theorem}[cf.~Theorem~13.2 of \cite{Titch}] \label{thm:main}
  Let $\epsilon > 0 $ be an arbitrary positive number.
  For a fixed non-negative integer $m$, the followings are equivalent~:
  \begin{enumerate}
  \item $\zeta^{(m)} (1 + it) = O(\log^\epsilon t)$.
  \item For each positive integer $k = 1, 2, \ldots$,
    \begin{equation*}
      \dfrac{1}{2\Delta} \int^{T+\Delta}_{T-\Delta} |\zeta^{(m)}(1+it)|^{2k} dt = O(\log^\epsilon T),
    \end{equation*}
    where $\Delta = 1/\log^{m+2} T$.
  \end{enumerate}
\end{theorem}

\begin{proof}
(1) $\Rightarrow$ (2) is obvious.

In order to show (2) $\Rightarrow$ (1), suppose that
$\zeta^{(m)}(1 + it)$ is not $O(\log^\epsilon t)$.
Then there exists $\lambda > 0$ and a sequence of complex numbers
$1+it_{\nu}$, such that $t_{\nu} \to \infty$ as $\nu \to \infty$,
and
\begin{equation*}
  |\zeta^{(m)}(1+it_{\nu})| > C \log^\lambda t_{\nu}
\end{equation*}
with some $C>0$.

On the other hand, we know
\begin{equation*}
  |\zeta^{(m+1)}(1 + it)| \leq E\log^{m+2} t
\end{equation*}
for sufficiently large $t$ and some constant $E$
(See the proof of Theorem~3.5 in \cite{Titch} and its remark).
Hence
\begin{align*}
  |\zeta^{(m)}(1+it) - \zeta^{(m)}(1+it_{\nu}) |
  &= \left| \int^t_{t_{\nu}} \zeta^{(m+1)}(1+iu) du \right| \\
  & \leq E \int^t_{t_{\nu}} |\log^{m+2} u| du \\
  &\leq A | t - t_{\nu} | \log^{m+2} t_{\nu} 
  < \frac{1}{2} C \log^\lambda t_{\nu},
\end{align*}
if $| t - t_{\nu} | \leq 1/\log^{m+2} t_{\nu} $ and $\nu$ is sufficiently large.
Hence
\begin{equation*}
  |\zeta^{(m)}(1+it)| > \frac{1}{2} C \log^\lambda t_{\nu}
  \qquad (| t - t_{\nu} | \leq 1/\log^{m+2} t_{\nu}).
\end{equation*}
Take $T = t_{\nu}$.
Then
\begin{align*}
  \frac{1}{2\Delta} \int^{T+\Delta}_{T-\Delta} |\zeta^{(m)}(1+it)|^{2k} dt
  &> \frac{1}{2\Delta} \int^{T+\Delta}_{T-\Delta} \left(\frac{1}{2} C \log^\lambda t_{\nu}\right)^{2k}dt \\
  &= \frac{1}{2\Delta}\left(\frac{1}{2} C\right)^{2k}  (\log^{2k\lambda} t_{\nu}) \cdot (2\Delta)
  = \left(\frac{1}{2} C\right)^{2k}  \log^{2k\lambda} t_{\nu}.
\end{align*}

This contradicts the assumption.
\end{proof}

\bigskip
{\small
  \begin{flushleft}
    Yoshihiro Koya\\
    \medskip
    Institute of Natural Science,\\
    Yokohama City University,\\
    22-2 Seto, Kanazawa-ku,\\
    Yokohama 236-0027, JAPAN\\
    \medskip
    \url{koya@yokohama-cu.ac.jp}
  \end{flushleft}}
\end{document}